\documentclass{article}

\usepackage{latexsym}
\usepackage{amsmath}
\usepackage{amsthm}
\usepackage{amssymb}
\usepackage{verbatim}
\usepackage{dsfont}

\newcommand{\mf}[1]{\mathfrak{#1}}

\newcommand{\mr}[1]{\mathrm{#1}}
\newcommand{\mb}[1]{\mathbb{#1}}

\newcommand{\eqdi}{\stackrel{d}{=}}

\newcommand{\norm}[1]{\| #1 \|}
\newcommand{\Norm}[1]{\bigl\| #1 \bigr\|}

\newcommand{\Sett}[1]{\bigl\{ #1 \bigr\}}

\newcommand{\sub}{\subseteq}

\newcommand{\ind}{\mathds{1}}

\newcommand{\oo}{\infty}

\newcommand{\RR}{\mb{R}}
\newcommand{\NN}{\mb{N}}

\newcommand{\EE}{\mb{E}\,}

\newcommand{\PP}{\mb{P}\,}

\newcommand{\XX}{\mf{N}}
\newcommand{\Xh}{\widehat{X}}

\newcommand{\al}{\alpha}
\newcommand{\be}{\beta}
               
\newcommand{\de}{\delta}               

\newcommand{\ve}{\varepsilon}

\newcommand{\la}{\lambda}

\newcommand{\vp}{\varphi}

\newcommand{\dist}{\mathrm{dist}\,}

\newtheorem{Counter}{!!!!!!!}

\newtheorem{Thm}[Counter]{Theorem}
\newtheorem{Pro}[Counter]{Proposition}
\newtheorem{Prop}[Counter]{Proposition}
\newtheorem{Cor}[Counter]{Corollary}
\newtheorem{Lem}[Counter]{Lemma}
\newtheorem{Rem}[Counter]{Remark}

\newcommand{\eav}{e^{\mr{av}}}

\begin{document}

\title{Free-Knot Spline Approximation of Stochastic Processes}
\author{
Jakob Creutzig\thanks{Fachbereich Mathematik, Technische Universit\"at
Darmstadt, Schlo\ss gartenstra\ss e 7, 64289 Darmstadt, Germany}%
\and 
Thomas M\"uller-Gronbach\thanks{Fakult\"at f\"ur Mathematik und Informatik,
FernUniversit\"at Hagen, L\"utzowstra\ss e 125, 58084 Hagen, Germany} 
\addtocounter{footnote}{-2}
\and 
Klaus Ritter\footnotemark 
\addtocounter{footnote}{1}
\and ~\\~\\ \sc \vspace{-17pt} \vspace{20pt} Dedicated to Henryk Wo\'zniakowski\\
\sc on the occasion of his 60th birthday
}

\date{December 6, 2006}

\maketitle

\begin{abstract}
We study optimal approximation of stochastic processes by polynomial
splines with free knots. The number of free knots is either a priori
fixed or may depend on the particular trajectory. For the $s$-fold
integrated Wiener process as well as for scalar diffusion processes
we determine the asymptotic behavior of the average $L_p$-distance
to the splines spaces, as the (expected) number $k$ of free knots
tends to infinity. 
\medskip

\end{abstract}

\section{Introduction}
Consider a stochastic process $X=(X(t))_{t\ge 0}$ with continuous
paths on a probability space $(\Omega,\mf{A},\PP)$. We study optimal
approximation of $X$ on the unit interval by polynomial splines with
free knots, which has first been treated in \cite{KP}.

For $k \in \NN$ and $r \in \NN_0$ we let $\Pi_r$ denote the set of
polynomials of degree at most $r$, and we consider the space
$\Phi_{k,r}$ of polynomial splines
\[
\vp = \sum_{j=1}^k \ind_{\left]t_{j-1}, t_j\right]} \cdot \pi_j,
\]
where $0 = t_0 < \ldots < t_k = 1$ and $\pi_1,\dots ,\pi_k \in
\Pi_r$. Furthermore, for $1\le p\le\oo$ and $1\le q <\oo$ we let
$\XX_{k,r}$ denote the class of measurable mappings
\[
\Xh: \Omega\to \Phi_{k,r}
\]
with $\Phi_{k,r}$ being considered as a subset of $L_p[0,1]$, and we
define
\[
e_{k,r} (X,L_p,q) = \inf\bigl\{
\bigl(\EE\|X-\Xh\|_{L_p[0,1]}^q\bigr)^{1/q}:
\Xh\in\XX_{k,r}\bigr\}.
\]
This quantity is an average $L_p$-distance from $X$ to
the space $\Phi_{k,r}$.

A natural extension of this methodology is
not to work with an a priori chosen number of free knots, but only
to control the average number of knots needed. This leads to the
definition $\Phi_r = \bigcup_{k=1}^\infty \Phi_{k,r}$ and to the
study of the class $\XX_{r}$ of measurable mappings
\[
\Xh:  \Omega \to \Phi_r
\]
with $\Phi_r$ being considered as a subset of $L_p[0,1]$. For a
spline approximation method $\Xh \in \XX_{r}$ we define
\[
\zeta(\Xh) = \EE ( \min \{ k \in \NN : \Xh(\cdot) \in \Phi_{k,r}\}),
\]
i.e., $\zeta(\Xh)-1$ is the expected number of free knots used by
$\Xh$. Subject to the bound $\zeta(\Xh)\le k$, the minimal achievable
error for approximation of $X$ in the class $\XX_{r}$ is given by
\[
\eav_{k,r} (X,L_p,q) = \inf\bigl\{
\bigl(\EE\|X-\Xh\|_{L_p[0,1]}^q\bigr)^{1/q}: \Xh\in\XX_{r},\
\zeta(\Xh) \le k\bigr\}.
\]
We shall study the asymptotics of the quantities $e_{k,r}$ and
$\eav_{k,r}$ as $k$ tends to infinity.

The spline spaces $\Phi_{k,r}$ form nonlinear manifolds that consist
of $k$-term linear combinations of functions of the form
$\ind_{\left]t,1\right]} \cdot \pi$ with $0 \leq t < 1$ and $\pi \in
\Pi_r$. Hence we are addressing a so-called nonlinear approximation
problem. While nonlinear approximation is extensively studied for
deterministic functions, see \cite {D} for a survey, much
less is known for stochastic processes, i.e., for random functions.
Here we refer to \cite{CA,CDGO}, where wavelet methods are analyzed,
and to \cite{KP}.
In the latter paper nonlinear approximation is related to
approximation based on partial information, as studied in
information-based complexity,
and spline approximation with free knots is analyzed as a
particular instance.

\section{Main Results}

For two sequences $(a_k)_{k
\in \NN}$ and $(b_k)_{k \in \NN}$ of positive real numbers we write
$a_k \approx b_k$ if $\lim_{k\to\oo}a_k / b_k =1$, and $a_k \gtrsim
b_k$ if $\liminf_{k\to\oo} a_k / b_k \ge 1$.  Additionally, $a_k
\asymp b_k$ means $c_1\le a_k/b_k \le c_2$ for all $k \in \NN$ and
some positive constants $c_i$.

Fix $s\in\NN_0$ and let $W^{(s)}$ denote an $s$-fold integrated
Wiener process. In \cite{KP}, the following result was proved.

\begin{Thm}\label{KP}
For $r\in\NN_0$ with $r\ge s$,
$$
  e_{k,r}(W^{(s)}, L_\oo, 1) \asymp
  \eav_{k,r}(W^{(s)}, L_\oo, 1) \asymp k^{-(s + 1/2)}.
$$
\end{Thm}

Our first result refines and extends this theorem. Consider the
stopping time
\[
\tau_{r,s,p}= \inf \Sett{t > 0 :
  \inf_{\pi \in \Pi_r}\|W^{(s)} - \pi\|_{L_p[0,t]} >  1},
\]
which yields the length of the maximal subinterval
$[0,\tau_{r,s,p}]$ that permits best approximation of $W^{(s)}$
from $\Pi_r$ with error at most one.
We have $0 < \EE \tau_{r,s,p} < \infty$, and we put
\[
\beta = s+1/2+1/p
\]
as well as
\[
c_{r,s,p} = (\EE \tau_{r,s,p})^{-\beta}
\]
and
\[
b_{s,p} = (s+1/2)^{s+1/2} \cdot p^{-1/p} \cdot \beta^{-\beta}.
\]

\begin{Thm}\label{thm:nonlinear-WP}
Let $r\in\NN_0$ with $r\ge s$ and $1 \leq q < \infty$. Then,
for $p=\oo$,
\begin{equation}\label{aus1}
\eav_{k,r}(W^{(s)}, L_\oo , q) \approx 
e_{k,r}(W^{(s)}, L_\oo , q) \approx   c_{r,s,\oo} \cdot k^{-(s+1/2)}.
\end{equation}
Furthermore, for $1 \leq p < \oo$,
\begin{equation}\label{aus2}
b_{s,p}
\cdot c_{r,s,p} \cdot k^{-(s+1/2)}
  \lesssim e_{k,r} (W^{(s)}, L_p, q)
  \lesssim
c_{r,s,p} \cdot k^{-(s+1/2)}
\end{equation}
and
\begin{equation}\label{aus3}
  \eav_{k,r}(W^{(s)}, L_p, q) \asymp k^{-(s+1/2)}.
\end{equation}
\end{Thm}

Note that the bounds provided by \eqref{aus1} and \eqref{aus2}
do not depend on the averaging parameter $q$. In
particular, asymptotic constants cannot explode for $q$ tending to
infinity.
Furthermore,
\[
\lim_{p \to \oo} b_{s,p} = 1
\]
for every $s \in \NN$,
but
\[
\lim_{s \to \oo} b_{s,p} = 0
\]
for every $1 \leq p < \oo$.
We conjecture that the upper bound in (i) is sharp.

We have an explicit construction of methods $\Xh^*_k \in \XX_{k,r}$
that achieve the upper bounds in \eqref{aus1} and \eqref{aus2},
i.e.,
\begin{equation}\label{explicit1}
\bigl(\EE\|W^{(s)} -
\widehat{X}^*_k\|_{L_p[0,1]}^q\bigr)^{1/q}\approx c_{r,s,p} \cdot
k^{-(s+1/2)},
\end{equation}
see \eqref{g0}.
Moreover, these methods a.s.~satisfy
\begin{equation}\label{explicit2}
\|W^{(s)} - \widehat{X}^*_k\|_{L_p[0,1]}\approx c_{r,s,p} \cdot
k^{-(s+1/2)}
\end{equation}
as well, while
\begin{equation}\label{explicit3}
\|W^{(s)} - \widehat{X}_k\|_{L_p[0,1]}\gtrsim b_{s,p} \cdot
c_{r,s,p} \cdot k^{-(s+1/2)}
\end{equation}
holds a.s.~for every sequence of approximations
$\widehat{X}_k\in\XX_{k,r}$. Here, $b_{s,\oo}=1$.

Our second result deals with approximation of a scalar diffusion
process given by the stochastic differential equation
\begin{equation}  \label{eq:sde-def}
\begin{aligned}
\phantom{,\qquad t\ge 0}
  dX(t) &= a(X(t))\, dt + b(X(t)) \, dW(t), \qquad t\ge 0, \\
  X(0) &= x_0.
\end{aligned}
\end{equation}
Here $x_0\in\RR$, and $W$ denotes a one-dimensional Wiener process.
Moreover, we assume that the functions $a,b:\RR\to\RR$ satisfy
\begin{itemize}
\item[(A1)] $a$ is Lipschitz continuous,
\item[(A2)] $b$ is differentiable with a bounded and Lipschitz
continuous derivative,
\item[(A3)] $b(x_0)\neq 0$.
\end{itemize}

\begin{Thm}\label{thm:nonlinear-diffusion}
Let $r\in\NN_0$, $1 \leq q < \oo$, and $1 \leq p \leq \oo$. Then
\[
e_{k,r}(X, L_p, q) \asymp \eav_{k,r}(X, L_p, q) \asymp k^{-1/2}
\]
holds for the strong solution $X$ of equation \eqref{eq:sde-def}.
\end{Thm}

For a diffusion process $X$ piecewise linear interpolation with free
knots is frequently used in connection with adaptive step-size
control. Theorem \ref{thm:nonlinear-diffusion} provides a lower
bound for the $L_p$-error of any such numerical algorithm, no matter
whether just Wiener increments or, e.g., arbitrary multiple
It\^{o}-integrals are used. Error estimates in \cite{HMGR,MG1} lead
to refined upper bounds in Theorem \ref{thm:nonlinear-diffusion} for
the case $1 \leq p <\oo$, as follows. Put
\[
\kappa(p_1,p_2)
= \bigl(\EE \|b\circ X\|_{L_{p_1}[0,1]}^{p_2}\bigr)^{1/p_2}
\]
for $1\le p_1,p_2 <\oo$. Furthermore, let $B$ denote a Brownian
bridge on
$[0,1]$ and define
\[
\eta(p) = \bigl(\EE\|B\|_{L_p[0,1]}^p\bigr)^{1/p}.
\]
Then
\[
e_{k,1} (X,L_p,p) \lesssim \eta(p)\cdot \kappa(2p/(p+2),p) \cdot k^{-1/2}
\]
and
\[
\eav_{k,1} (X,L_p,p) \lesssim \eta(p)\cdot\kappa(2p/(p+2),2p/(p+2))
\cdot k^{-1/2}.
\]
We add that these upper bounds are achieved by numerical
algorithms with adaptive step-size control for the
Wiener increments.

In the case $p=\oo$ it is interesting to compare the results on
free-knot spline approximation with average $k$-widths of $X$.
The latter quantities are defined by
\[
d_k (X,L_p,q) = \inf_\Phi \Bigl( \EE \Bigl( \inf_{\varphi \in \Phi}
\|X - \varphi\|^q_{L_p [0,1]} \Bigr) \Bigr)^{1/q},
\]
where the infimum is taken over all linear subspaces $\Phi \subseteq
L_p[0,1]$ of dimension at most $k$.
For $X=W^{(s)}$ as well as in the diffusion case we have
\[
d_k (X,L_\oo,q) \asymp k^{-(s+1/2)},
\]
see \cite{C,M1,M2,M3} and \cite{DMGR}.
Almost optimal linear subspaces are not known
explicitly, since the proof of the upper bound for
$d_k (X,L_\oo,q)$ is non-constructive.
We add that in the case of an $s$-fold integrated Wiener process
piecewise polynomial interpolation of $W^{(s)}$ at
equidistant knots $i/k$ only yields errors of order
$(\ln k)^{1/2} \cdot k^{-(s+1/2)}$,
see \cite{R} for results and references.
Similarly, in the diffusion case, methods
$\Xh_k\in\XX_r$ that are only based on
pointwise evaluation of $W$ and satisfy $\zeta(\Xh_k)\le k$ can at
most achieve errors of order $(\ln k)^{1/2} \cdot k^{-1/2}$,
see \cite{MG2}.

\section{Approximation of Deterministic Functions}

Let $r \in \NN_0$ and $1 \leq p \leq \oo$ be fixed. We   introduce
error measures, which allow to determine suitable free knots for
spline approximation. For $f \in C\left[0, \oo\right[$ and $0 \leq u
< v$ we put
\[
  \delta_{[u,v]}(f) = \inf_{\pi \in \Pi_r}  \|f - \pi\|_{L_p[u,v]}.
\]
Furthermore, for $\ve > 0$, we put $\tau_{0, \ve}(f) = 0$, and we
define
\[
  \tau_{j, \ve}(f) = \inf\{ t > \tau_{j-1, \ve}(f) : \,
  \delta_{[\tau_{j-1, \ve}(f),t]}(f) > \ve\}
\]
for $j \geq 1$. Here $\inf \emptyset = \oo$, as usual.
Put $I_j(f) = \{\ve > 0 :\, \tau_{j,\ve}(f) < \oo\}$.

\begin{Lem}\label{l1}
Let $j \in\NN$.
\begin{itemize}
\item[(i)]
If $\ve \in I_j(f)$ then
\[
\delta_{[\tau_{j-1,\ve}(f),\tau_{j,\ve}(f)]} (f) = \ve.
\]
\item[(ii)]
The set $I_j(f)$ is an interval, and the mapping $\ve \mapsto
\tau_{j,\ve}(f)$ is strictly increasing and right-continuous
on $I_j(f)$. Furthermore,
$\tau_{j,\ve}(f) > \tau_{j-1,\ve}(f)$ if $\ve \in I_{j-1}(f)$,
and
$\lim_{\ve\to\oo} \tau_{j,\ve}(f) = \oo$.
\item[(iii)]
If $v \mapsto \delta_{[u,v]}(f)$ is strictly increasing for
every $u\ge 0$,
then $\ve \mapsto\tau_{j,\ve}(f)$ is continuous on $I_j(f)$.
% vgl. Revusz, Yor p. 8, 179
\end{itemize}
\end{Lem}

\begin{proof}
First we show that the mapping $(u,v) \mapsto \delta_{[u,v]} (f)$ is
continuous.
Put $J_1 = [u/2,u+(v-u)/3]$ as well as $J_2 = [v-(v-u)/3,2v]$.
Moreover, let $\pi^\alpha (t) = \sum_{i=0}^r \alpha_i \cdot t^i$ for
$\alpha \in \RR^{r+1}$, and define a norm on $\RR^{r+1}$ by
\[
\|\alpha\| = \|\pi^\alpha\|_{L_p[u+(v-u)/3,v-(v-u)/3]}.
\]
If $(x,y)\in J_1\times J_2$ and
\[
\|f - \pi^\alpha\|_{L_p[x,y]} = \delta_{[x,y]}(f)
\]
then
\[
\|\alpha\| \le \|\pi^\alpha\|_{L_p[x,y]} \le
 \delta_{[u/2,2v]}(f) + \|f\|_{L_p[u/2,2v]}.
\]
Hence there exists a compact set $K \subseteq \RR^{r+1}$ such that
\[
\delta_{[x,y]} (f) = \inf_{\alpha \in K} \|f -
\pi^\alpha\|_{L_p[x,y]}
\]
for every $(x,y) \in J_1 \times J_2$.
Since $(x,y,\alpha) \mapsto \|f - \pi^\alpha\|_{L_p[x,y]}$ defines a
continuous mapping on $J_1 \times J_2 \times K$, we conclude that
$(x,y) \mapsto \inf_{\alpha \in K} \|f - \pi^\alpha\|_{L_p[x,y]}$ is
continuous, too, on $J_1 \times J_2$.

Continuity and monotonicity of $v \mapsto \delta_{[u,v]}(f)$
immediately imply (i).

The monotonicity stated in (ii) will be verified inductively. Let $0
< \ve_1 < \ve_2 $ with $\ve_2 \in I_j(f)$, and suppose that
$\tau_{j-1,\ve_1}(f) \leq \tau_{j-1,\ve_2}(f)$. Note that the latter
holds true by definition for $j=1$. {}From (i) we get
\[
\delta_{[\tau_{j-1,\ve_1}(f),\tau_{j,\ve_2}(f)]} (f) \geq
\delta_{[\tau_{j-1,\ve_2}(f),\tau_{j,\ve_2}(f)]} (f) = \ve_2.
\]
This implies $\tau_{j,\ve_1}(f) \leq \tau_{j,\ve_2}(f)$, and (i)
excludes equality to hold here.

Since $\delta_{[u,v]}(f)\le \|f\|_{L_p[u,v]}$,
the mappings $\ve \mapsto \tau_{j,\ve}(f)$ are unbounded
and $\tau_{j,\ve}(f) > \tau_{j-1,\ve}(f)$ if $\ve \in I_{j-1}(f)$.

For the proof of the continuity properties stated in (ii) and (iii)
we also proceed inductively, and we use (i) and the monotonicity
from (ii).
Consider a sequence $(\ve_n)_{n \in \NN}$
   in $I_j(f)$, which converges monotonically to $\ve\in I_j(f)$, and put
$t= \lim_{n \to \oo} \tau_{j,\ve_n}(f)$. Assume that
$\lim_{n\to\oo}\tau_{j-1,\ve_n}(f)=\tau_{j-1,\ve}(f)$, which
obviously holds true for $j=1$. Continuity of $(u,v) \mapsto
\delta_{[u,v]}(f)$ and (i) imply
$\delta_{[\tau_{j-1,\ve}(f),t]}(f) = \ve$, so that $t \leq
\tau_{j,\ve}(f)$. For a decreasing sequence $(\ve_n)_{n \in \NN}$
we also have $\tau_{j,\ve}(f) \leq t$. For an increasing sequence
$(\ve_n)_{n \in \NN}$ we use the strict monotonicity of $v \mapsto
\delta_{[u,v]}(f)$ to derive $t= \tau_{j,\ve}(f)$.
\end{proof}

Let $F$ denote the class of functions $f\in C\left[0,\oo\right[$
that satisfy
\begin{equation}\label{finite}
\tau_{j,\ve}(f) < \oo
\end{equation}
for every $j\in\NN$ and $\ve >0$ as well as
\begin{equation}\label{g21}
\lim_{\ve \to 0} \tau_{j,\ve}(f) = 0
\end{equation}
for every $j \in \NN$.

Let $k\in\NN$. We now present an almost optimal spline approximation
method of degree $r$  with $k-1$ free knots for functions $f\in F$.
Put
\[
\gamma_k(f) = \inf \{ \ve > 0 : \tau_{k,\ve}(f)  \geq 1 \}
\]
and note that \eqref{g21} together with Lemma \ref{l1}.(ii)
implies $\gamma_k(f)\in\left]0,\oo\right[$.
Let
\[
\tau_j = \tau_{j,\gamma_k(f)}(f)
\]
for $j=0,\ldots,k$ and define
\begin{equation}\label{g0}
  \varphi_k^*(f) =
\sum_{j=1}^k \ind_{\left]\tau_{j-1}, \tau_j\right]}
  \cdot
  \mr{argmin}_{\pi \in \Pi_r}  \|f - \pi\|_{L_p[\tau_{j-1}, \tau_j]}.
\end{equation}
Note that Lemma \ref{l1} guarantees
\begin{equation}\label{g7}
\|f - \varphi_k^*(f)\|_{L_p[\tau_{j-1}, \tau_j]} = \gamma_k(f)
\end{equation}
for $j=1,\dots,k$ and
\begin{equation}\label{g5}
\tau_k  \geq 1.
\end{equation}

The spline $\varphi_k^*(f)|_{[0,1]} \in \Phi_{k,r}$ enjoys the following
optimality properties.

\begin{Prop}\label{psi-psi^*-Wiener}
Let $k \in \NN$ and $f\in F$.
\begin{itemize}
\item[(i)]
For $1 \leq p \leq \infty$,
$$
  \|f - \varphi_k^*(f)\|_{L_p[0,1]}
\leq  k^{1/p} \cdot \gamma_k(f).
$$
\item[(ii)]
For $p=\oo$ and every $\varphi\in \Phi_{k,r}$,
\[
\|f - \varphi\|_{L_{\oo}[0,1]} \ge \gamma_{k}(f).
\]
\item[(iii)]
For $1\le p < \oo$, every $\varphi \in \Phi_{k,r}$,
and every $m \in \NN$ with $m > k$,
\[
\|f - \varphi\|_{L_p[0,1]} \ge (m-k)^{1/p} \cdot \gamma_{m}(f).
\]
\end{itemize}
\end{Prop}

\begin{proof}
For $p < \oo$,
\[
  \|f - \varphi_k^*(f)\|_{L_p[0,1]}^p
 \leq
  \sum_{j=1}^k \|f - \varphi_k^*(f)\|_{L_p[\tau_{j-1},\tau_j]}^p
= k \cdot (\gamma_k(f))^p
\]
follows from \eqref{g7} and \eqref{g5}. For $p=\oo$, (i) is verified
analogously.

Consider a polynomial spline
$\varphi\in\Phi_{k,r}$ and let $0=t_0  < \ldots < t_k = 1$
denote the corresponding knots. Furthermore, let
$\rho\in\left]0,1\right[$.
For the proof of (ii) we put
\[
\sigma_j = \tau_{j,\rho\cdot\gamma_k(f)}(f).
\]
for $j=0,\ldots,k$. Then $\sigma_k < 1$, which implies
\[
[\sigma_{j - 1}, \sigma_{j}]\sub [t_{j-1}, t_{j}]
\]
for some $j\in\{1,\ldots,k\}$. Consequently, by Lemma \ref{l1},
\[
  \|f - \varphi\|_{L_{\oo}[0,1]}
  \ge
  \|f - \varphi\|_{L_{\oo}[\sigma_{j-1},
     \sigma_{j}]}
  \ge \inf_{\pi \in \Pi_r}
 \|f - \pi\|_{L_{\oo}[\sigma_{j-1},
     \sigma_{j}]}
 = \rho\cdot\gamma_{k}(f).
\]

For the proof of (iii) we define
\[
\sigma_\ell = \tau_{\ell,\rho\cdot\gamma_{m}(f)}(f)
\]
for $\ell=0,\dots,m$. Then $\sigma_{m} < 1$, which implies
\[
[\sigma_{\ell_i - 1}, \sigma_{\ell_i}] \sub [t_{j_i-1},
t_{j_i}]
\]
for some indices $1\le j_1 \le \ldots \le j_{m-k} \le k$ and
$1\le \ell_1 < \ldots < \ell_{m-k} \le m$. Hence, by Lemma \ref{l1},
\[
\|f - \varphi\|_{L_p[0,1]}^{p}
  \ge \sum_{i=1}^{m-k}
  \inf_{\pi \in \Pi_r} \|f - \pi\|^p_%
{L_p[\sigma_{\ell_i-1},\sigma_{\ell_i}]}
  = (m-k) \cdot \rho^p \cdot (\gamma_{m}(f))^p.
\]
for $1\le p < \oo$. Letting $\rho$ tend to one completes the proof.
\end{proof}

\section{Approximation of Integrated Wiener Processes}

Let $W$ denote a Wiener process and consider the $s$-fold integrated
Wiener processes $W^{(s)}$ defined by $W^{(0)} = W$ and
\[
  W^{(s)}(t) = \int_0^t W^{(s-1)}(u) \, du
\]
for $t \geq 0$ and $s\in\NN$. We briefly discuss some properties of
$W^{(s)}$, that will be important in the sequel.

The scaling property of the Wiener process implies that for every
$\rho > 0$ the process $(\rho^{-(s+1/2)} \cdot W^{(s)}(\rho\cdot
t))_{t \geq 0}$ is an $s$-fold integrated Wiener process, too. This
fact will be called the scaling property of $W^{(s)}$.

While $W^{(s)}$ has no longer independent increments for $s \geq 1$,
the influence of the past is very explicit. For $z >0$ we define
${_zW}^{(s)}$ inductively by
\[
  {_zW}^{(0)}(t) = W(t + z) - W(z)
\]
and
\[
  {{_zW}^{(s)}(t)} = {\int_0^{t}} {_zW}^{(s-1)}(u)\, du.
\]
Then it is easy to check that
\begin{equation}\label{intWP-increments}
W^{(s)}(t + z) = \sum_{i=0}^s \frac{t^i}{i!} W^{(s-i)}(z) +
{_zW}^{(s)}(t).
\end{equation}

Consider the filtration generated by $W$, which coincides with the
filtration generated by $W^{(s)}$, and let $\tau$ denote a stopping
time with $\PP(\tau < \oo)=1$. Then the strong Markov property of
$W$ implies that the process
\[
{_{\tau}W}^{(s)} = ({_{\tau}W}^{(s)}(t))_{t \ge 0}
\]
is an $s$-fold integrated Wiener process, too. Moreover, the
processes ${_{\tau}W}^{(s)}$ and $(\ind_{[0, \tau]}(t) \cdot
W(t))_{t \ge 0}$ are independent, and consequently, the processes
${_{\tau}W}^{(s)}$ and $(\ind_{[0, \tau]}(t) \cdot W^{(s)}(t))_{t
\ge 0}$ are independent as well. These facts will be called the
strong Markov property of $W^{(s)}$.

Fix $s\in\NN_0$. In the sequel we assume that $r \geq s$.
For any fixed $\ve > 0$ we consider the sequence
of stopping times $\tau_{j,\ve}(W^{(s)})$,
which turn out to be
finite a.s.~and therefore are strictly increasing, see Lemma~\ref{l1}.
Moreover, for $j\in\NN$,
we define
\[
\xi_{j,\ve} = \tau_{j,\ve}(W^{(s)}) - \tau_{j-1,\ve}(W^{(s)}).
\]
These random variables
yield the lengths of
consecutive maximal subintervals that permit best approximation from
the space $\Pi_r$ with error at most $\ve$.
Recall that $F \subseteq C\left[0,\oo\right[$ is defined via
properties \eqref{finite} and \eqref{g21} and that
$\beta = s+1/2+1/p$.

\begin{Lem}\label{l2}
The $s$-fold integrated Wiener process $W^{(s)}$ satisfies
\[
\PP(W^{(s)}\in F) = 1.
\]
For every $\ve > 0$ the random variables $\xi_{j,\ve}$
form an i.i.d.~sequence with
\[
   \xi_{1, \ve} \eqdi \ve^{1/\beta}
   \cdot \xi_{1, 1}\qquad\text{and}\qquad \EE ( \xi_{1,1}) < \oo.
\]
\end{Lem}

\begin{proof}
We claim that
\begin{equation}\label{moment}
\EE(\tau_{j,\ve}(W^{(s)})) < \oo
\end{equation}
for every $j\in\NN$.

For the case $j=1$ let $Z = \delta_{[0,1]} (W^{(s)})$ and note that
\[
\delta_{[0,t]} (W^{(s)}) \eqdi t^{\beta} \cdot Z
\]
follows for $t > 0$ from the scaling property of $W^{(s)}$.
Hence we have
\begin{equation}\label{g9}
\PP(\tau_{1,\ve}(W^{(s)}) < t) =
\PP(\delta_{[0,t]} (W^{(s)}) > \ve) =
\PP(Z > \ve \cdot t^{-\beta}),
\end{equation}
which, in particular, yields
\begin{equation}\label{g4}
\tau_{1,\ve}(W^{(s)}) \eqdi \ve^{1/\beta} \cdot
\tau_{1,1}(W^{(s)}).
\end{equation}
According to Corollary \ref{app:smdev-del-WP},
there exists a constant $c>0$ such that
\[
\PP(Z \leq \eta) \leq \exp(-c \cdot \eta^{-1/(s+1/2)})
\]
holds for every $\eta \in \left]0,1\right]$.
We conclude that
\[
\PP(\tau_{1,1}(W^{s)}) > t) \leq \exp(-c \cdot t)
\]
if $t \geq 1$, which implies $\EE (\tau_{1,1}(W^{(s)})) < \oo$.

Next, let $j \geq 2$,  put  $\tau = \tau_{j-1,\ve}(W^{(s)})$
and $\tau^\prime = \tau_{j,\ve}(W^{(s)})$, and assume
that $\EE(\tau)<\oo$.
{}From the representation \eqref{intWP-increments}
and the fact that $r \ge s$ we derive
\[
  \de_{[\tau, t]}(W^{(s)}) = \de_{[0, t - \tau]}
    ({_{\tau} W^{(s)}}),
\]
and hence it follows that
\begin{equation}\label{ind}
  \tau^\prime = \tau
  + \tau_{1, \ve}({_{\tau} W^{(s)}}).
\end{equation}
We have $\EE(\tau_{1,\ve}({_{\tau}W^{(s)}})) < \oo$,
since ${_{\tau} W^{(s)}}$ is an $s$-fold integrated Wiener process
again, and consequently
$\EE(\tau^\prime)<\oo$.

We turn to the properties of the sequence $\xi_{j,\ve}$. Due to
\eqref{g4} and \eqref{ind} we have
\[
\xi_{j,\ve} = \tau_{1, \ve}({_{\tau} W^{(s)}}) \eqdi
\tau_{1, \ve}(W^{(s)}) \eqdi \ve^{1/\beta}\cdot \xi_{1,1}.
\]
Furthermore, $\xi_{j,\ve}$ and
$(\ind_{[0, \tau]}(t) \cdot W^{(s)}(t))_{t \ge 0}$
are independent because of the strong Markov property of $W^{(s)}$,
and therefore
$\xi_{j,\ve}$ and $(\xi_{1,\ve},\ldots,\xi_{j-1,\ve})$ are
independent as well.

It remains to show that the trajectories of $W^{(s)}$ a.s.~satisfy
\eqref{g21}. By the properties of the sequence $\xi_{j,\ve}$ we have
\begin{equation}\label{g8}
\tau_{j, \ve}(W^{(s)})\eqdi \ve^{1/\beta}\cdot \tau_{j,1}(W^{(s)}).
\end{equation}
Observing \eqref{moment} we conclude that
\begin{align*}
\PP ( \lim_{\ve \to 0} \tau_{j,\ve}(W^{(s)}) \geq t )
&=
\lim_{\ve \to 0} \PP (\tau_{j,\ve}(W^{(s)}) \geq t)\\
&=
\lim_{\ve \to 0} \PP (\tau_{j,1}(W^{(s)}) \geq t / \ve^{1/\beta}) = 0
\end{align*}
for every $t>0$, which completes the proof.
\end{proof}

Because of Lemma \ref{l2},
Proposition \ref{psi-psi^*-Wiener}
yields sharp upper and lower bounds for
the error of spline approximation of $W^{(s)}$ in terms of the
random variable
\[
V_k = \gamma_k(W^{(s)}).
\]

\begin{Rem}\label{equal}
{\rm
Note that $W^{(s)}$
a.s.~satisfies $W^{(s)}|_{[u,v]} \not\in \Pi_r$
for all $0 \leq u < v$.
Assume that $p < \oo$. Then
$v \mapsto \delta_{[u,v]}(W^{(s)})
$ is a.s.~strictly increasing for all $u \geq 0$.
We use Lemma \ref{l1}.(iii) and Lemma \ref{l2} to conclude that,
with probability one, $V_k$ is the unique solution of
\[
\tau_{k,V_k}(W^{(s)}) = 1.
\]

Consequently, we a.s.~have equality
in Proposition \ref{psi-psi^*-Wiener}.(i)
for $1\le p < \infty$, too.
Note that with positive probability solutions $\ve$ of the equation
$\tau_{k,\ve}(W^{(s)}) = 1$ fail to exist in the case $p=\infty$.
}
\end{Rem}

To complete the analysis of spline approximation methods
we study the asymptotic behavior of the sequence $V_k$.

% Man koennte uns beim folgenden Lemma und in Appendix A
% Dilettantismus vorwerfen. Etwa: Noch nie von
% uniform integrability gehoert, was?

\begin{Lem}\label{l4}
For every $1 \leq q < \oo$,
\[
\left(\EE V_k^q \right)^{1/q} \approx (k \cdot \EE (\xi_{1,1}))^{-\beta}.
\]
Furthermore, with  probability one,
\[
V_k\approx (k \cdot \EE (\xi_{1,1}))^{-\beta}.
\]
\end{Lem}

\begin{proof}
   Put
\[
S_k = 1/k \cdot \sum_{j=1}^k \xi_{j, 1}
\]
and use \eqref{g8} to obtain
\begin{equation}\label{g13}
\PP ( V_k \leq \ve)  =
\PP ( \tau_{k,\ve}(W^{(s)}) \geq 1)
=
\PP ( k^{-\beta} \cdot S_k^{-\beta} \leq \ve).
\end{equation}
Therefore
\[
\EE (V_k^q)  =
k^{-\beta q} \cdot  \EE ( S_k^{-\beta q}),
\]
and for the first statement it remains to show that
\[
\EE ( S_k^{-\beta q}) \approx (\EE
(\xi_{1,1}))^{-\beta q}.
\]

The latter fact follows from Proposition \ref{momentsofmeans},
if we can verify
that $\xi_{1,1}$ has a proper lower tail behavior \eqref{gt}.
To this end we use \eqref{g9} and the large deviation
estimate \eqref{ldev-Wiener} to obtain
\begin{align*}
\PP ( \xi_{1,1} \leq \eta) & =
\PP ( \delta_{[0,1]} (W^{(s)}) \geq \eta^{-\beta} )\\
& \leq
\PP ( \|W^{(s)}\|_{L_p[0,1]} \geq \eta^{-\beta} ) \\
& \leq \exp( - c \cdot \eta^{-2\beta})
\end{align*}
with some constant $c>0$ for all $\eta \leq 1$.

In order to prove the second statement, put
\[
S_k^* = (k \cdot \sigma^2)^{-1/2} \cdot\sum_{j=1}^k (\xi_{j,1} -
\mu),
\]
where $\mu = \EE (\xi_{1,1})$ and $\sigma^2$ denotes the variance
of $\xi_{1,1}$. Let $\rho > 1$. Then
\[
\PP (V_k > \rho \cdot (k \cdot \mu)^{-\beta})
=
\PP (S_k < \rho^{-1/\beta} \cdot \mu)
=
\PP (S_k^* < k^{1/2} \cdot \widetilde{\rho})
\]
with
\[
\widetilde{\rho} = (\rho^{-1/\beta}-1)/\sigma \cdot \mu < 0,
\]
due to \eqref{g13}.
We apply a local version of the central limit theorem,
which holds for i.i.d.~sequences with a finite third moment,
see~\cite[Thm.~V.14]{Pe},
to obtain
\begin{align*}
&\PP (V_k > \rho \cdot (k \cdot \mu)^{-\beta}) \\
&\qquad \leq
c_1 \cdot k^{-1/2} \cdot (1+k^{1/2} \cdot |\widetilde{\rho}|)^{-3}
+ (2\pi)^{-1/2} \cdot \int_{-\oo}^{k^{1/2} \cdot \widetilde{\rho}}
\exp(-u^2/2)\, du\\
&\qquad\leq
c_2 \cdot k^{-2}
\end{align*}
with constants $c_i > 0$.
For every $\rho < 1$ we get
\begin{equation}\label{clt}
\PP (V_k < \rho \cdot (k \cdot \mu)^{-\beta})
\leq c_2 \cdot k^{-2}
\end{equation}
in the same way.
It remains to apply the Borel-Cantelli Lemma.
\end{proof}

Observe Remark \ref{equal} and use Lemma \ref{l4} to obtain
\eqref{explicit1} and \eqref{explicit2} for the methods $\Xh^*_k =
\varphi_k^*(W^{(s)})$. Clearly, \eqref{explicit1} implies the 
upper bounds in \eqref{aus1}, \eqref{aus2}, and \eqref{aus3}.

For the proof of the lower bound in \eqref{aus2}
we apply Proposition \ref{psi-psi^*-Wiener}.(iii) with
\[
m = \lfloor  \beta / (s+1/2) \cdot k \rfloor.
\]
Then we have
\begin{align*}
&(m-k)^{1/p} \cdot \left(\EE V_m^q \right)^{1/q}\\
&\qquad
\approx
k^{-(s+1/2)} \cdot p^{-1/p}
\cdot \beta^{-\beta} \cdot
(s+1/2)^{s+1/2} \cdot
(\EE (\xi_{1,1}))^{-\beta},
\end{align*}
as claimed. Using the second statement in Lemma \ref{l4},
the lower bound in \eqref{explicit3} is shown in the same
way.

It remains to prove the lower bounds for $\eav_{k,r}(W^{(s)},L_p,q)$
in Theorem \ref{thm:nonlinear-WP}.

\begin{proof}[Proof of the lower bound in \eqref{aus1}]
Let $k \in \NN$ and consider $\Xh_k \in \XX_r$ such that 
$\zeta (\Xh_k) \leq k$, i.e.,
\begin{equation}\label{eq:newav:1}
  \sum_{\ell = 1}^{\oo} \ell \cdot \PP(B_\ell) \le k
\end{equation}
for $B_\ell = \Sett{\widehat{X}(\cdot) \in \Phi_{\ell, r}
\setminus \Phi_{\ell-1,r}}$, 
where $\Phi_{0,r} = \emptyset$. 
By Proposition \ref{psi-psi^*-Wiener}.(ii),
\[
\EE \Norm{W^{(s)} - \Xh_k}_{L_{\oo}[0,1]}^q 
  \ge 
  \sum_{\ell=1}^{\oo} \EE (\ind_{B_\ell} \cdot V_\ell^q).
\]
For $\varrho \in \left]0,1\right[$, $\mu = \EE (\xi_{1,1})$,
and $L \in \NN$ we define 
$$
A_\ell = \Sett{V_\ell > \rho \cdot (\ell \cdot \mu)^{-\be}}, 
$$
and
\[
C_L = \bigcup_{\ell=1}^L B_\ell.
\]
Since $\gamma_\ell(f) \geq \gamma_{\ell+1}(f)$ for $f \in F$,
we obtain
\begin{align*} 
&\EE \Norm{W^{(s)} - \Xh_k}_{L_{\oo}[0,1]}^q \\
&\qquad \ge 
\sum_{\ell=1}^L \EE(\ind_{B_\ell} \cdot V_L^q) 
+ \sum_{\ell = L+1}^\oo \EE(\ind_{B_\ell} \cdot V_\ell^q) \\  
&\qquad \ge \sum_{\ell=1}^L \EE(\ind_{B_\ell \cap A_L} \cdot V_L^q) 
+ \sum_{\ell=L+1}^\oo 
\EE(\ind_{B_\ell \cap A_\ell} \cdot V_\ell^q)\\ 
&\qquad \ge \rho^q \mu^{-\beta q} \cdot
\Bigl(L^{- \be q} \cdot  \PP(C_L \cap A_L) + \sum_{l = L+1}^\oo
\ell^{-\be q} \cdot \PP(B_\ell \cap A_\ell)\Bigr).
\end{align*}
{}From \eqref{clt}
we infer that $\PP(A_\ell^c) \le c_1 \cdot \ell^{-2}$ 
with a constant $c_1>0$. 
Hence there exists a constant $c_2 >0$ such that
\[
\Gamma (L) =
L^{- \be q}\cdot \PP(C_L) 
  + \sum_{\ell = L+1}^\oo \ell^{-\be q} \cdot \PP(B_\ell)
  - c_2 \cdot L^{-\be q - 1}
\]
satisfies
\begin{equation}\label{eq:newav:4} 
\rho^{-q}  \mu^{\beta q} \cdot
\EE \Norm{W^{(s)} - \Xh_k}_{L_{\oo}[0,1]}^q 
\ge \Gamma(L)
\end{equation}
for every $L \in \NN$.

Put $\al = (1 + 2\be q)/(2 + 2\be q)$, and take 
$L(k) \in [k^{\al}-1, k^{\al}]$.  We claim that there 
exists a constant $c_3 > 0$  such that 
\begin{equation}\label{eq:newav:goal} 
  k^{\be q}\cdot \Gamma (L(k)) \ge 
  \left(1 - k^{-(1 - \al)\be q}\right)^{1 + \be q} 
  - c_3 \cdot k^{-1/2}.
\end{equation}

First, assume that
$\PP(C_L) \ge k^{-(1 - \al)\be q}$.  Then 
\begin{align*}
  k^{\be q} \cdot \Gamma (L(k)) &\ge 
  k^{\be q}   \cdot \Bigl(k^{-\al \be q} \cdot
\PP(C_L) - c_2 \cdot (k^\alpha-1)^{- \be q-1}\Bigr) \\
  &\ge  1 - c_3 \cdot k^{-1/2}
\end{align*}
with a constant $c_3>0$.
Next, assume $\PP(C_L) < k^{-(1 - \al)\be q}$ and use \eqref{eq:newav:1} 
to derive
\begin{align*}
1- k^{-(1 - \al)\be q}
& < 
\sum_{\ell=L+1}^\oo \PP(B_\ell)\\
&=
\sum_{\ell=L+1}^\oo (\ell \cdot\PP(B_\ell))^{\be q/(1+\be
q)} \cdot  (\ell^{-\be q} \cdot \PP (B_\ell))^{1/(1+\be q)}
 \\
& \leq
\biggl(\, \sum_{\ell=L+1}^\oo \ell \cdot \PP (B_\ell) 
\biggr)^{\be q/(1+\be q)}
\cdot
\biggl(\, \sum_{\ell=L+1}^\oo \ell^{-\be q}
 \cdot \PP (B_\ell) \biggr)^{1/(1+\be q)}\\
&\leq 
k^{\be q/(1+\be q)} \cdot
\biggl(\, \sum_{\ell=L+1}^\oo \ell^{-\be q}
 \cdot \PP (B_\ell) \biggr)^{1/(1+\be q)}.
\end{align*}
Consequently,
\begin{align*}
k^{\be q}\cdot \Gamma (L(k)) 
&\ge 
k^{\be q}\cdot
\biggl(\,
  \sum_{\ell = L+1}^\oo \ell^{-\be q} \cdot \PP(B_\ell)
  - c_2 \cdot (k^\alpha-1)^{- \be q-1}
\biggr)\\
&\ge
  \left(1 - k^{-(1 - \al)\be q}\right)^{1 + \be q} 
- c_3 \cdot k^{-1/2},
\end{align*}
    which completes the proof of \eqref{eq:newav:goal}.
By \eqref{eq:newav:4} and \eqref{eq:newav:goal}, 
\[
\EE \Norm{W^{(s)} - \Xh_k}_{L_{\oo}[0,1]}^q 
\gtrsim
\rho^{q}  \mu^{-\beta q} \cdot k^{-\be q}
\]
for every $\rho \in \left]0,1\right[$.
\end{proof}

Finally, for the proof of the 
lower bound in \eqref{aus3} it suffuces to 
establish the lower bound claimed for $\eav_{k,r}(W^{(s)},L_1,1)$.
For further use, we shall prove a more general result.

\begin{Lem}\label{l99}
For every $s \in \NN$ there exists a constant $c>0$
with the following property.
For every $\Xh\in\mf{N}_r$, every $A \in\mf{A}$ with $\PP (A) \geq 4/5$,
and every $t \in \left]0,1\right]$ we have
\[
\EE \left( \ind_{A} \cdot \|W^{(s)} - \Xh \|_{L_1[0,t]}\right)
\geq c \cdot t^{s+3/2} \cdot k^{-(s+1/2)}.
\]
\end{Lem}

\begin{proof}
Because of the scaling property of $W^{(s)}$ it
suffices to study the particular case $t=1$.
Put
\[
B = \{ \Xh \in \Phi_{2k,r} \},
\]
and observe that $\PP (B) \geq 1/2$ follows from $\xi(\Xh)\le k$.
Due to Lemma \ref{l2} and Proposition \ref{psi-psi^*-Wiener}.(iii),
\[
\ind_B \cdot \| W^{(s)} - \Xh \|_{L_1[0,1]} \geq
\ind_B \cdot 2k \cdot V_{4k}.
\]

Put $\mu =  \EE (\xi_{1,1})$, choose
$0 < c < ( 2 \mu)^{-\beta}$, and define
\[
D_k = \{V_k \geq c \cdot k^{-\beta} \}.
\]
By \eqref{g13} we obtain
\[
\PP (D_k)  = \PP ( S_k < c^{-1/\beta} ) \geq
\PP ( S_k < 2 \mu ).
\]
Hence
\[
\lim_{k\to\infty} \PP(D_k) =1
\]
due to the law of large numbers, and consequently
$\PP (B \cap D_k) \geq 2/5$
if $k$ is sufficiently large, say $k \geq k_0$.
We conclude that
\[
\ind_{A \cap B\cap D_{4k}}
\cdot \| W^{(s)} - \Xh \|_{L_1[0,1]} \geq
\ind_{A \cap B\cap D_{4k}} \cdot c \cdot 2^{1-2\beta} \cdot k^{-(s+1/2)}
\]
and
$\PP(A \cap B \cap D_{4k}) \geq 1/5$ if $4k \geq k_0$.
Take expectations to complete the proof.
\end{proof}

Lemma \ref{l99} with
$A = \Omega$ and $t= 1$
yields the lower bound in \eqref{aus3}

\section{Approximation of Diffusion Processes}

Let $X$ denote the solution of the stochastic differential equation
\eqref{eq:sde-def} with initial value $x_0$, and recall that the
drift coefficient $a$ and the diffusion coefficient $b$ are supposed
to satisfy conditions (A1)--(A3). In the following we use $c$ to
denote unspecified positive constants, which may only depend on
$x_0$, $a$, $b$ and the averaging parameter
$1 \leq q < \oo$.

Note that
\begin{equation}\label{eq:X-bounded}
  \EE \norm{X}_{L_{\oo}[0,1]}^q < \oo
\end{equation}
and
\begin{equation}\label{eq:X-Holder}
  \EE \Bigl(\, \sup_{t \in [s_1,s_2]} |X(t) - X(s_1)|^q \Bigr)
\le c \cdot (s_2-s_1)^{q/2}
\end{equation}
   for all $1\le q <\oo$ and  $0\le s_1\le s_2 \le 1$, see \cite[p.\ 138]{KlPl}.

\subsection{Upper Bounds}

In order to establish upper bounds, it suffices to consider the case
of $p=\oo$ and $r=0$, i.e., nonlinear approximation in supremum norm
with piecewise constant splines.

\begin{comment}
For convenience we extend $X$ to a stochastic process on
$\left[0,\oo\right[$ by setting $X(t) := X(1) + W(t) - W(1)$ for $t
> 1$. Furthermore, we define
\[
H(f)(t) := f(t) - x_0 - \int_0^{t\wedge 1} a(f(s)) \, \mr{d} s
\]
for $f \in C\left[0,\oo\right[$ and $t \geq 0$, and hereby we
dissect $X$ into its martingale part $M := H(X)$ and $Y := X -
H(X)$.
%$$
%  Y(t) := x_0 + \int_0^{t \wedge 1} a(X(s)) \, \mr{d} s
%$$
%and
%$$
%M(t) :=
%\int_0^{t \wedge 1} b(X(s)) \, \mr{d} W(s) + W(t) - W(t \wedge 1),
%$$
%where $t \geq 0$.
\end{comment}

We dissect $X$ into its martingale part
\[
M(t) = \int_0^t b(X(s)) \, dW(s)
\]
and
\[
Y(t) =  x_0 + \int_0^t a(X(s)) \, ds.
\]

\begin{Lem}\label{Lem:drift}
For all $1 \leq q < \oo$ and $k \in \NN$, there exists an
approximation $\widehat{Y} \in \XX_{k,0}$ such that
$$
  \left(\EE \norm{Y - \widehat{Y}}_{L_\oo[0,1]}^q \right)^{1/q}
\leq c \cdot k^{-1}.
$$
\end{Lem}

\begin{proof}
Put $\|g\|_{\mr{Lip}} = \sup_{0 \leq s< t \leq 1} |g(t)-g(s)|/|t-s|$
  for $g:[0,1]\to\RR$, and define
$$
\widehat{Y} = \sum_{j=1}^{k} \ind_{](j-1)/k, j/k]} \cdot Y((j-1)/k).
$$
     By (A1)             and \eqref{eq:X-bounded},
$$
  \EE \norm{Y - \widehat{Y}}_{L_\oo[0,1]}^q
\le \EE \norm{Y}_{\mr{Lip}}^q \cdot k^{-q}
  \le  c\cdot
\bigl(1 + \EE \norm{X}_{L_\oo[0,1]}^q\bigr) \cdot k^{-q}
   \le c\cdot k^{-q}.
$$
\end{proof}

\begin{Lem}
For all $1 \leq q < \oo$ and $k \in \NN$, there exists an
approximation $\widehat{M} \in \XX_{k,0}$ such that
$$
  \Bigl(\EE \norm{M - \widehat{M}}_{L_\oo[0,1]}^q \Bigr)^{1/q}
\leq c \cdot k^{-1/2}.
$$
\end{Lem}

\begin{proof}
  Let
$$
  \Xh = \sum_{j=1}^{k} \ind_{\left](j-1)/k, j/k\right]}
\cdot X((j-1)/k).
$$
Clearly, by \eqref{eq:X-Holder},
$$
  \Bigl(\EE \|X - \Xh\|_{L_2[0,1]}^q\Bigr)^{1/q}
\leq c \cdot k^{-1/2}.
$$

     Define
$$
  R(t) = \int_0^t b(\Xh(s)) \, dW_s.
$$
By the Burkholder-Davis-Gundy inequality and (A2),
\begin{align}\label{g76}
\notag
   \Bigl( \EE \| M - R \|_{L_\oo[0,1]}^q\Bigr)^{1/q}
   &\le c \cdot \biggl( \EE \Bigl(\int_0^1
\notag
   (b(X(s)) - b(\Xh (s)))^2 \, ds \Bigr)^{q/2} \biggr)^{1/q}  \\
   &\le c \cdot
\notag
  \Bigl(\EE \|X - \Xh\|_{L_2[0,1]}^q\Bigr)^{1/q} \\
&\leq c \cdot k^{-1/2}.
\end{align}
   Note that
$$
  R = \widehat{R} + V,
$$
where
$$
  \widehat{R} = \sum_{j = 1}^{k}
  \ind_{\left](j-1)/k, j/k\right]} \cdot R((j-1)/k)
$$
and
$$
  V = \sum_{j=1}^{k}  \ind_{\left](j-1)/k, j/k\right]}  \cdot
            b(X((j-1)/k)) \cdot (W - W((j-1)/k)).
$$
According to Theorem \ref{thm:nonlinear-WP}, there exists an
approximation $\widehat{W} \in \XX_{k,0}$ such that
$$
  \Bigl( \EE \|W - \widehat{W}\|_{L_{\oo}[0,1]}^{2q} \Bigr)^{1/(2q)}
  \leq c \cdot k^{-1/2}.
$$
   Using $\widehat{W}$ we define
$\widehat{V} \in \XX_{2k,0}$ by
\[
  \widehat{V} = \sum_{j=1}^{k}  \ind_{\left](j-1)/k, j/k\right]}  \cdot
            b(X((j-1)/k)) \cdot (\widehat{W} - W((j-1)/k)).
\]
    Clearly,
$$
  \|V - \widehat{V}\|_{L_{\oo}[0,1]} \le \|b(X)\|_{L_{\oo}[0,1]} \cdot
  \|W - \widehat{W}\|_{L_{\oo}[0,1]}.
$$
Observing \eqref{eq:X-bounded} and (A2), we conclude that
\begin{align}\label{g77}
\notag
  \Bigl( \EE \|V - \widehat{V}\|^q_{L_{\oo}[0,1]}\Bigr)^{1/q}
 &\le \Bigl( \EE \|b(X)\|^{2q}_{L_{\oo}[0,1]}\Bigr)^{1/(2q)} \cdot
 \Bigl( \EE \|W - \widehat{W}\|^{2q}_{L_{\oo}[0,1]}\Bigr)^{1/(2q)} \\
&\leq c \cdot k^{-1/2}.
\end{align}

We finally define $\widehat{M} \in \XX_{2k,0}$ by $\widehat{M} =
\widehat{R} + \widehat{V}$.
Since
\[
M-\widehat{M}= (M-R)+(V-\widehat{V}),
\]
it remains to apply the estimates \eqref{g76} and \eqref{g77}
to complete the proof.
\end{proof}

The preceding two lemma imply $e_{k,0}(X, L_\oo, q) \leq c
\cdot k^{-1/2}$.

\subsection{Lower Bounds}

For establishing lower bounds it suffices to study the case
$p=q=1$.
Moreover, we assume without loss of generality that $b(x_0) >0$.

Choose
$\eta >0$ as well as a function $b_0: \RR \to \RR$ such that
\begin{itemize}
\item[(a)] $b_0$ is differentiable with a bounded and Lipschitz
continuous derivative,
\item[(b)]  $\inf_{x \in \RR} b_0(x) \ge b(x_0)/2$,
\item[(c)]  $b_0 = b$ on the interval
$[x_0 - \eta, x_0 + \eta]$.
\end{itemize}

    We will use a Lamperti transform based on the
space-transformation
$$
  g(x) = \int_{x_0}^{x} \frac{1}{b_0(u)}\, du.
$$
   Note that $g'= 1/b_0$ and $g'' = - b_0'/b_0^2$,
and define $H_1,H_2: C[0,\oo[\to C[0,\oo[$ by
$$
  H_1(f)(t) = \int_0^t
\bigl(g'a + g''/2 \cdot b^2\bigr)(f(s))\, ds
$$
and
$$
  H_2(f)(t) = g(f(t)).
$$
   Put $H = H_2 - H_1$. Then by the It\^o formula,
$$
  H(X) (t) = \int_0^t \frac{b(X(s))}{b_0(X(s))} \, dW(s).
$$

The idea of the proof is as follows. We show that any good spline
approximation of $X$ leads to a good spline approximation of $H(X)$.
However, since with a high  probability, $X$ stays within
$[x_0 - \eta, x_0 + \eta]$ for some short (but nonrandom) period of
time, approximation of $H(X)$ is not easier than
approximation of $W$, modulo constants.

First, we consider approximation of $H_1(X)$.

\begin{Lem}\label{l6}
For every $k \in \NN$ there exists an approximation $\Xh_1 \in
\XX_{k,0}$ such that
$$
\EE \norm{H_1(X) - \Xh_1}_{L_1[0,1]} \leq c \cdot k^{-1}.
$$
\end{Lem}

\begin{proof}
Observe that $\bigl|g'a + g''/2 \cdot b^2\bigr|(x) \leq c \cdot
(1+x^2)$, and proceed as in the Proof of Lemma \ref{Lem:drift}.
\end{proof}

Next, we relate approximation of $X$ to approximation of $H_2(X)$.

\begin{Lem}\label{HX-to-X}
For every approximation $\Xh \in \XX_r$ there exists an approximation
$\Xh_2 \in \XX_r$ such that
\[
\zeta(\Xh_2) \le 2 \cdot \zeta(\Xh)
\]
and
\[
\EE \|H_2(X) - \Xh_2\|_{L_1[0,1]} \le c\cdot \bigl( \EE \|X -
\Xh\|_{L_1[0,1]} + 1/\zeta(\Xh) \bigr).
\]
\end{Lem}

\begin{proof}
For a fixed $\omega \in \Omega$ let $\Xh (\omega) $ be given by
$$
\Xh(\omega) = \sum_{j=1}^{k} \ind_{]t_{j-1}, t_j]} \cdot \pi_j.
$$
We refine the corresponding partition to a partition
$0=\widetilde{t}_0 < \ldots < \widetilde{t}_{\widetilde{k}} = 1$
that contains all the points $i/\ell$, where $\ell = \zeta(\Xh)$.
Furthermore, we define the polynomials $\widetilde{\pi}_j \in \Pi_r$
by
$$
  \Xh(\omega) = \sum_{j=1}^{\widetilde{k}}
\ind_{\left]\widetilde{t}_{j-1}, \widetilde{t}_j\right]} \cdot
\widetilde{\pi}_j.
$$
Put $f = X(\omega)$ and define
\[
\Xh_2(\omega)= \sum_{j=1}^{\widetilde{k}}
\ind_{\left]\widetilde{t}_{j-1}, \widetilde{t}_j \right]} \cdot q_j
\]
with polynomials
$$
  q_j =
g(f(\widetilde{t}_{j-1})) + g'(f(\widetilde{t}_{j-1})) \cdot
(\widetilde{\pi}_j - f(\widetilde{t}_{j-1})) \in \Pi_r.
$$

Let $\widehat{f}_2 = \Xh_2(\omega)$.
If $t \in \left]\widetilde{t}_{j-1}, \widetilde{t}_j\right] \subseteq
\left](i-1)/\ell,i/\ell\right]$, then
\begin{align*}
  &|H_2 (f)(t) - \widehat{f}_2(t)| \\
&\qquad = \bigl| g(f(t)) - g(f(\widetilde{t}_{j-1}))
 - g'(f(\widetilde{t}_{j-1})) \cdot (\widetilde{\pi}_j(t) -
  f(\widetilde{t}_{j-1}))\bigr| \\
  &\qquad \le
\left| g(f(t)) - g(f(\widetilde{t}_{j-1}))
 - g'(f(\widetilde{t}_{j-1})) \cdot
(f(t) - f(\widetilde{t}_{j-1}))\right| \\
&\qquad \phantom{\le\ }
 + \left|g'(f(\widetilde{t}_{j-1}))\right| \cdot
|f(t) - \widetilde{\pi}_j(t)| \\
&\qquad \leq c \cdot \left( |f(t) - f(\widetilde{t}_{j-1})|^2 +
|f(t) - \widetilde{\pi}_j(t)| \right) \\
&\qquad \leq c \cdot \Bigl(\, \sup_{t \in
\left](i-1)/\ell,i/\ell\right]} |f(t) - f((i-1)/\ell)|^2 + |f(t) -
\widetilde{\pi}_j(t)| \Bigr).
\end{align*}
Consequently, we may invoke \eqref{eq:X-Holder} to derive
\[
\EE \|H_2(X) - \Xh_2\|_{L_1[0,1]} \leq c \cdot \bigl(
1/\zeta(\Xh) + \EE \|X - \Xh\|_{L_1[0,1]} \bigr).
\]
Moreover, $\zeta(\Xh_2) \leq 2 \cdot \zeta(\Xh)$.
\end{proof}

Finally, we establish a lower bound for approximation of $H(X)$.

\begin{Lem}\label{l88}
For every approximation $\Xh \in \XX_r$,
\[
\EE \| H(X) - \Xh\|_{L_1[0,1]} \geq c \cdot (\zeta(\Xh))^{-1/2}.
\]
\end{Lem}

\begin{proof}
Choose $t_0 \in \left]0,1\right]$ such that
\[
A = \Bigl\{ \sup_{t \in [0,t_0]} |X(t) - x_0| \le \eta \Bigr\}
\]
satisfies $\PP(A) \geq 4/5$. Observe that
\[
\ind_{A} \cdot \|H(X) - \Xh\|_{L_1[0,1]} \geq \ind_{A} \cdot \|W -
\Xh\|_{L_1[0,t_0]},
\]
and apply Lemma \ref{l99} for $s= 0$.
\end{proof}

\begin{proof}[Proof of the lower bound in Theorem
\ref{thm:nonlinear-diffusion}] Consider any approximation $\Xh \in
\XX_r$ with $k-1 < \zeta(\Xh) \leq k$, and choose $\Xh_1$ and
$\Xh_2$ according to Lemma \ref{l6} and Lemma \ref{HX-to-X},
respectively. Then
\begin{align*}
&\EE \| H(X) - (\Xh_2 - \Xh_1)\|_{L_1[0,1]}\\
&\qquad\leq \EE \| H_2(X) - \Xh_2\|_{L_1[0,1]} +
\EE \| H_1(X) - \Xh_1\|_{L_1[0,1]}\\
&\qquad\leq c \cdot \bigl( \EE \|X - \Xh\|_{L_1[0,1]} +
(\zeta(\Xh))^{-1} + k^{-1} \big)\\
&\qquad\leq c \cdot \bigl( \EE \|X - \Xh\|_{L_1[0,1]} + k^{-1}
\bigr).
\end{align*}
On the other hand, $\zeta(\Xh_2 - \Xh_1) \leq \zeta(\Xh_2) + k
\leq 3 \cdot k$, so that
\[
\EE \| H(X) - (\Xh_2 - \Xh_1)\|_{L_1[0,1]} \geq c \cdot k^{-1/2}
\]
follows from Lemma \ref{l88}. We conclude that
\[
\EE \|X - \Xh\|_{L_1[0,1]} \geq c\cdot k^{-1/2},
\]
as claimed.
\end{proof}

\appendix

\section{Convergence of Negative Moments of Means}

Let $(\xi_i)_{i \in \NN}$ be an i.i.d.~sequence of
random variables such that $\xi_1 >0$ a.s.~and $\EE (\xi_1) < \oo$.
Put
\[
S_k = 1/k \cdot \sum_{i=1}^{k} \xi_i.
\]

\begin{Pro}\label{momentsofmeans}
For every  $\al > 0$,
\[
\liminf_{k \to \oo} \EE (S_k^{-\al}) \ge (\EE (\xi_1))^{-\al}.
\]
If
\begin{equation}\label{gt}
\phantom{\qquad\quad v \in \left]0,v_0\right],}
\PP(\xi_1 < v) \le c \cdot v^{\rho},
\qquad\quad v \in \left]0,v_0\right],
\end{equation}
for some constants $c, \rho, v_0 > 0$, then
$$
 \lim_{k \to \oo} \EE (S_k^{-\al}) = (\EE (\xi_1))^{-\al}.
$$
\end{Pro}

\begin{proof}
Put $\mu = \EE( \xi_1)$ and define
\[
g_k(v) = \al \cdot v^{-(\al +1)} \cdot \PP(S_k < v).
\]
Thanks to the weak law of large numbers, $\PP(S_k < v)$
tends to $\ind_{\left]\mu, \oo\right[}(v)$ for every
every $v\neq \mu$. Hence, by Lebesgue's theorem,
\begin{equation}\label{eq:pf:Pro1:1}
\lim_{k \to \oo} \int_{\mu/2}^{\oo}
g_k(v) \, dv = \mu^{-\al}\; .
\end{equation}

Since
\[
\EE (S_k^{-\al}) = \int_0^{\oo} \PP(S_k^{-\al} > u)\, du
               = \int_0^{\oo} g_k(v) \, dv
\]
the asymptotic lower bound for $\EE (S_k^{-\al})$ follows from
\eqref{eq:pf:Pro1:1}.

Given \eqref{gt},
we may assume without loss of generality that
$c \cdot v_0^{\rho} < 1$.
We first consider the case $\xi_1\le 1$ a.s.,
and we put
\[
A_k = \int_{x_0/k}^{\mu/2} g_k(v) \, dv\qquad\text{and}\qquad
  B_k = \int_{0}^{x_0/k} g_k(v) \, dv.
\]

For $v_0/k \le v \le \mu/2$ we use Hoeffding's inequality to obtain
\[
g_k(v) \le  \al \cdot v^{-(\al +1)} \cdot \PP(|S_k - \mu| > \mu/2)
\le  \al \cdot (k/v_0)^{\al+1}\cdot 2\exp(-  k/2 \cdot \mu^2),
\]
which implies
\[
\lim_{k\to\infty} A_k =0.
\]
On the other hand, if $\rho k > \al$, then
\begin{align*}
  B_k &= k^{\al+1} \cdot \al \cdot \int_{0}^{v_0} v^{- (\al +1)} \cdot
    \PP\Bigl(\sum_{i=1}^k \xi_i < v \Bigr) \, dv\\
&\le k^{\al+1} \cdot \al \cdot \int_{0}^{v_0} v^{- (\al +1)} \cdot
    (\PP(\xi_1 < v))^k\, dv\\
&\le k^{\al+1} \cdot \al \cdot c^k \cdot
\int_{0}^{v_0} v^{\rho k - (\al +1)}
\, dv\\
&= k^{\al+1} \cdot \al \cdot (\rho k - \al)^{-1} \cdot c^k \cdot
    v_0^{\rho k - \al},
\end{align*}
and therefore
\[
\lim_{k\to\infty} B_k =0.
\]
In view of \eqref{eq:pf:Pro1:1}
we have thus proved the proposition in the case of bounded variables $\xi_i$.

In the general case put $\xi_{i,N} = \min\{N, \xi_i\}$ as well as
$S_{k,N} = 1/k \cdot  \sum_{i=1}^k \xi_{i,N}$,
and apply the result for bounded variables to obtain
\[
  \limsup_{k \to \oo}\EE(S_k^{-\al})
  \le \inf_{N \in \NN} \limsup_{k \to \oo} \EE(S_{k,N}^{-\al})
  = \inf_{N \in \NN} (\EE \xi_{1,N})^{-\al}
  = (\EE \xi_1)^{- \al}
\]
by the monotone convergence theorem.
\end{proof}

\section{Small Deviations of $W^{(s)}$ from $\Pi_r$}\label{ab}

Let $X$ denote a centered Gaussian random variable with values in
a normed space $(E,\|\cdot\|)$, and consider a finite-dimensional
linear subspace $\Pi\subset E$. We are interested in the small
deviation behavior of
\[
d(X,\Pi) = \inf_{\pi\in\Pi} \|X-\pi\|.
\]

Obviously,
\begin{equation}\label{g2}
\PP(\|X\| \le \ve) \leq \PP( d(X, \Pi) \le \ve)
\end{equation}
for every $\ve>0$.
We establish an upper bound for
$\PP( d(X, \Pi) \le \ve)$
that involves large deviations of $X$, too.

\begin{Pro}\label{p1}
If $\dim(\Pi)=r$ then
$$
  \PP( d(X, \Pi) \le \ve) \le
(4 \la/ \ve)^{r} \cdot \PP(\|X\| \le 2 \ve) +
  \PP(\|X\| \ge \la - \ve)
$$
for all $\la \ge \ve >0$.
\end{Pro}

\begin{proof}
Put $B_\delta(x) = \{y\in E:\, \|y-x\|\le\delta\}$ for $x\in E$ and
$\delta>0$, and consider the sets $A=\Pi\cap B_\la(0)$ and $B=B_\ve(0)$.
Then
\[
\{d(X,\Pi)\le\ve\} \subset \{X\in A+B\}\cap\{\|X\|\ge\la-\ve\},
\]
and therefore it suffices to prove
\begin{equation}\label{genug}
\PP(X\in A+B) \leq (4 \la/ \ve)^{r} \cdot \PP(\|X\|\leq 2\ve).
\end{equation}

Since $1/\la\cdot A\subset \Pi\cap B_1(0)$, the $\ve$-covering number
of $A$ is not larger than $(4\la/\ve)^r$, see \cite[Eqn. (1.1.10)]{CS}.
Hence
\[
A\subset\bigcup_{i=1}^n B_\ve(x_i)
\]
for some $x_1,\ldots,x_n\in E$ with $n\le (4\la/\ve)^r$, and
consequently,
\[
A+B\subset \bigcup_{i=1}^n B_{2\ve}(x_i).
\]
Due to Anderson's inequality we have
\[
\PP(X\in B_{2\ve}(x_i)) \le \PP(X\in B_{2\ve}(0)),
\]
which implies \eqref{genug}.
\end{proof}

Now, we turn to the specific case of $X=(W^{(s)}(t))_{t\in[0,1]}$
and $E=L_p[0,1]$, and we consider the subspace $\Pi=\Pi_r$ of
polynomials of degree at most $r$.

According to the large deviation principle for the
$s$-fold integrated Wiener process,
\begin{equation}\label{ldev-Wiener}
- \log \PP(\|W^{(s)}\|_{L_p[0,1]} > t) \asymp t^{2}
\end{equation}
as $t$ tends to infinity, see, e.g.,  \cite{DZ}. Furthermore,
the small ball probabilities satisfy
 \begin{equation}\label{smdev-Wiener}
- \log \PP(\|W^{(s)}\|_{L_p[0,1]} \le \ve) \asymp
\ve^{-1/(s+1/2)}
\end{equation}
as $\ve$ tends to zero, see, e.g., \cite{Li} and \cite{Lif}.

\begin{Cor}\label{app:smdev-del-WP}
For all $r,s\in\NN_0$,
\[
  - \log \PP(d(W^{(s)}, \Pi_r) \le \ve) \
\asymp  \ve^{-1/(s+1/2)}
\]
as $\ve$ tends to zero.
\end{Cor}

\begin{proof}
Use \eqref{g2} and \eqref{smdev-Wiener} to obtain the upper bound.
For the lower bound employ Proposition \ref{p1} with $\la = \ve^{-1}$,
and note that
\[
-\log \PP(\|W^{(s)}\|_{L_p[0,1]} \ge 1/\ve - \ve) \preceq
-\log \PP(\|W^{(s)}\|_{L_p[0,1]} \le 2\ve)
\]
as $\ve$ tends to zero, due to \eqref{ldev-Wiener}
and \eqref{smdev-Wiener}.
\end{proof}

\section*{Acknowledgments}
The authors are grateful to Mikhail Lifshits for helpful discussions.
In particular, he pointed out to us the approach in Appendix \ref{ab}.

\end{document}